\documentclass[12pt]{amsart}
\usepackage{amssymb}
\newcommand{\C}{\mathcal{C}}
\newcommand{\tF}{\tilde{F}}
\newcommand{\tE}{\tilde{E}}
\newcommand{\tg}{\tilde{g}}
\newcommand{\tG}{\tilde{G}}
\newcommand{\tC}{\tilde{\mathcal{C}}}
\newcommand{\tP}{\tilde{P}}

\newtheorem{Def}{Definition}
\newtheorem{Lem}{Lemma}
\newtheorem{Prop}{Proposition}
\newtheorem{Thm}{Theorem}
\newtheorem{Cor}{Corollary}
\newtheorem{Rem}{Remark}
\newenvironment{Pf}{ Proof.}{\(\square\)}

\title[On Asanov's Finsleroid-Finsler metrics...]{On Asanov's Finsleroid-Finsler metrics as the solutions of a conformal rigidity problem}
\author{Cs. Vincze}
\address{Inst. of Math., Univ. of Debrecen \\
H-4010 Debrecen, P.O.Box 12 \\
Hungary}
\email{csvincze@science.unideb.hu}
\begin{document}
\begin{abstract}
Finsleroid-Finsler metrics form an important class of singular (y-local) Finsler metrics. They were introduced by G. S. Asanov \cite{A2} in 2006. As the special case of the general construction Asanov produced singular (y - local) examples of  Landsberg spaces of dimension at least three that are not of Berwald type. The existence of regular (y - global) Landsberg metrics that are not of Berwald type is an open problem up to this day; for a detailed exposition of the so-called unicorn problem in Finsler geometry see D. Bao \cite{B}. 

In this paper we are going to characterize the Finsleroid-Finsler metrics as the solutions of a conformal rigidity problem. We are looking for (non-Riemannian) Finsler metrics admitting a (non-homothetic) conformal change such that the mixed curvature tensor of the Berwald connection contracted by the derivatives of the logarithmic scale function is invariant. We prove that the solutions of class at least $\mathcal{C}^2$ on the complement of the zero section are conformal to Finsleroid-Finsler metrics. 
\end{abstract}
\maketitle

\section*{Introduction}
In the paper we prove the following theorems.

\vspace{0.5cm}
\noindent
{\bf{Theorem A}} Let $M$ be a manifold of dimension $n\geq 3$ equipped with the Finslerian metric function $F$ and consider a function $\alpha \colon M \to \mathbb{R}$ satisfying the regularity condition $d_p\alpha\neq 0$ at a point $p\in M$. 
If the sectional curvature of the indicatrix hypersurface at $p$ is positive and 
$$\tilde{P}_{ijk}^l=P^l_{ijk}$$
under the conformal change $\tilde{F}=e^{\alpha}F$, where $P^l_{ijk}$ is the mixed curvature tensor of the Berwald connection, then $F$ is a locally Riemannian metric function.

\vspace{0.5cm}
\noindent
{\bf{Theorem B}} Let $M$ be a manifold of dimension $n\geq 3$ equipped with the Finslerian metric function $F$ and consider a function $\alpha \colon M \to \mathbb{R}$ satisfying the regularity condition $d_p\alpha\neq 0$ at a point $p\in M$.
If the sectional curvature of the indicatrix hypersurface at $p$ is positive and 
$$\tilde{P}_{ijk}^l\alpha_l=P^l_{ijk}\alpha_l$$
under the conformal change $\tilde{F}=e^{\alpha}F$, where $\alpha_l$ denotes the partial derivatives of $\alpha$ (depending only on the position\emph), then $F$ is locally conformal to a Finsleroid-Finsler metric or it is a locally Riemannian metric function. 

\vspace{0.5cm}
Looking for results like Theorem A and Theorem B was originally motivated by the so-called Matsumoto's problem \cite{M2} in 2001: are there conformally related (non-Riemannian) Berwald manifolds? The problem is closely related to the intrinsic characterization of Wagner manifolds \cite{V6}, see also \cite{V10}. Wagner manifolds form a special class of generalized Berwald manifolds admitting a linear connection on the base manifold such that the parallel transports preserve the Finslerian norm of tangent vectors. Especially, the compatible linear connection of a Wagner manifold is semi-symmetric with a special (exact or at least closed) one-form in the usual decomposition of its torsion. Wagner manifolds can be also defined as conformally Berwald Finsler manifolds due to M. Hashiguchi and Y. Ychijyo \cite{HY}, see also \cite{V1}, \cite{V2} and \cite{V9}. The logarithm of the scale function (the logarithmic scale function) corresponds to the potential of the one-form in the torsion of the compatible linear connection up to a constant proportional term. It is clear that the scale function of the conformal change of a Finslerian metric function to a Berwaldian one is uniquely determined (up to a constant homothetic term) if and only if the conformality between two (non-Riemannian) Berwald manifolds must be trivial (homothetic). The first attempt to solve Matsumoto's problem was given in \cite{V3}, where the author investigated the consequences of the conformal invariance of the mixed curvature tensor of the Berwald connection (see Theorem A). It is a natural generalization of the original problem\footnote{Bochner's technic and the theory of geometric vector fields in the tangent spaces of a Finsler manifold give another way to solve the generalized Matsumoto's problem in \cite{V7}.}. In what follows we generalize the basic results of \cite{V3} to prove Theorem B which gives the conformal characterization of Finsleroid-Finsler metrics \cite{A2}. The cronology of the basic steps:
\begin{itemize}
\item 1998 - the central symmetric version of the Finsleroid-Finsler metric in G. S. Asanov \cite{A1}.
\item 2003 - the non-symmetric version of the Finsleroid-Finsler metric in \cite{V3} as \emph{Asanov-type Finslerian metric functions}. They satisfy differential equation (\ref{diffeqversion3dricatti}) of the generalized Matsumoto's problem (conformal invariance of the mixed curvature tensor of the Berwald connection) for the Finslerian energy along special directions in the tangent spaces; see also \cite{V9}.
\item 2006 - Asanov's necessary and sufficient conditions for (non-symmetric) Finsleroid-Finsler metrics to be of Landsberg but not of Berwald type;
\item 2016 - non-symmetric Finsleroid-Finsler metrics as the solution of a conformal rigidity problem (the invariance of the {\bf{contracted}} mixed curvature tensor of the Berwald connection), see Theorem B and section 6.2. (the converse of the theorem);
\end{itemize}
for a detailed exposition of the unicorn problem in Finsler geometry see D. Bao \cite{B}. In what follows we give a characterization of Finsleroid-Finsler metrics as the singular non-Riemannian solutions of the conformal rigidity problem $\tilde{P}_{ijk}^l\alpha_l=P^l_{ijk}\alpha_l$, where $P_{ijk}^l$'s are the components of the mixed curvature tensor of the Berwald connection and $\alpha_l$'s are the partial derivatives of the "logarithmic" scale function $\alpha$ depending only on the position. The basic steps of the proof are Theorem 1, Theorem 2, the solution of a Ricatty-type diffrential equation (sections 4.2, 4.4, 4.5 and 4.6) and Theorem 6.

\section*{Acknowledgement}
The paper was motivated by the oral communication with Professor David Bao at the 50th Symposium on Finsler Geometry (21-25. Oct. 2015, Hiroshima, Japan). I would like to thank him for paying my attention to some correspondences between Asanov's Unicorn metrics and Finsler metrics satisfying conformal rigidity properties. I am very grateful for his human and professional encouragement.

The work is supported by the University of Debrecen's internal research project RH/885/2013.

\section{Notations and terminology}

Let $M$ be a manifold with local coordinates $u^1, \ldots, u^n.$ The induced coordinate system of the tangent manifold $TM$ consists of the functions
$$x^1, \ldots, x^n\ \ \textrm{and}\ \ y^1, \ldots, y^n,$$
where $x$'s refer to the coordinates of the base point and $y$'s denote the coordinates of the directions:
$$v\in T_pM \ \ \textrm{can be written as}\ \ v=y^{i}(v)\frac{\partial}{\partial u^i}_{x(v)},\ \ \textrm{where}\ \ x(v)=p.$$
\subsection{Finsler metrics} A Finsler metric is a continuous function $F\colon TM\to \mathbb{R}$ satisfying the following conditions:
\begin{itemize}
\item $\displaystyle{F}$ is smooth on the complement of the zero section (regularity),
\item $\displaystyle{F(tv)=tF(v)}$ for all $\displaystyle{t> 0}$ (positive homogenity),
\item the Hessian
$$g_{ij}=\frac{\partial^2 E}{\partial y^i \partial y^j}$$
of the Finslerian energy function $\displaystyle{E=(1/2)F^2}$ is positive definite at all nonzero elements $\displaystyle{v\in T_pM}$ (strong convexity). It is called the Riemann-Finsler metric of the Finsler manifold.
\end{itemize}
In what follows we summerize some basic notations and facts we need to prove our theorems. As a general reference of Finsler geometry and the forthcoming list of quantities see \cite{B} and \cite{BSC}:
\begin{itemize}
\item $\displaystyle{F}$ is a Finsler metric function, $\displaystyle{E:=\frac{1}{2}F^2}$ is the energy function,
\item $\displaystyle{l_i=\frac{\partial F}{\partial y^i}}$, 
\item $\displaystyle{g_{ij}=\frac{\partial^2 E}{\partial y^i \partial y^j}}$ is the Riemann-Finsler metric and its inverse $\displaystyle{g^{ij}=(g_{ij})^{-1}},$
\item $\displaystyle{C=y^l\frac{\partial}{\partial y^l}}$ is the Liouville vector field, 
\item $\displaystyle{C_{ijk}=\frac{\partial g_{ij}}{\partial y^k}}$
is the so-called first Cartan tensor\footnote{The Cartan tensor quantities are often defined as
$$A_{ijk}=\frac{F}{2}\frac{\partial g_{ij}}{\partial y^k}$$
see e.g. \cite{B} and \cite{BSC}. The symbol $\mathcal{C}_{ijk}$ follows \cite{G1} and \cite{G2}.}, $\mathcal{C}^l_{ij}=g^{lk}\mathcal{C}_{ijk}$. The first Cartan tensor is totally symmetric and $\displaystyle{y^k\C_{ijk}=0}$. It is also known that
\begin{equation}
\label{eq:1}
\frac{\partial g^{lm}}{\partial y^i}=-2g^{mk}\C_{ik}^l=-2\C_{i}^{lm}
\end{equation}
and, consequently,
\begin{equation}
\label{eq:2}
\frac{\partial \C_{jk}^l}{\partial y^i}-\frac{\partial \C_{ik}^l}{\partial y^j}=2\left(\C^l_{jm}\C_{ik}^m-\C^l_{im}\C_{jk}^m\right), \ \ \textrm{where}\ \ Q_{ijk}^l=\C^l_{jm}\C_{ik}^m-\C^l_{im}\C_{jk}^m
\end{equation}
is the $vv$-curvature of the Cartan connection. 
\end{itemize}
\subsection{Geodesic spray coefficients:}
\begin{equation}
\label{eq:3}
G^l=\frac{1}{2}g^{lm}\left(y^k\frac{\partial Fl_m}{\partial x^k}-F\frac{\partial F}{\partial x^m}\right),\ \ \textrm{i.e.}\ \ G^l=\frac{1}{2}g^{lm}\left(y^k\frac{\partial^2 E}{\partial y^m\partial x^k}-\frac{\partial E}{\partial x^m}\right).
\end{equation}

\subsection{Horizontal sections:}

$$\frac{\delta }{\delta x^i}=\frac{\partial }{\partial x^i}-G^l_i\frac{\partial }{\partial y^l},\ \ \textrm{where}\ \ G_i^l=\frac{\partial G^l}{\partial y^i}.$$

\subsection{The second Cartan tensor:} (Landsberg tensor\footnote{Following subsection 1.6. the Landsberg tensor is often defined as
$$\stackrel{\ .}{A}_{ijk}=-\frac{1}{2}Fl_i\frac{\partial G_{i}^l}{\partial y^j};$$
see e.g. \cite{B} and \cite{BSC}. The symbol $P_{ij}^l$ follows \cite{H2} and \cite{M1}.})
$$P^l_{ij}=\frac{1}{2} g^{lm}\left( \frac{\delta g_{jm}}{\delta x^i}-G^k_{ij}g_{km}-G^k_{im} g_{jk}\right)=\frac{1}{2} g^{lm}\left( \frac{\partial g_{jm}}{\partial x^i}-2G^k_i\C_{jkm}-G^k_{ij}g_{km}-G^k_{im} g_{jk}\right),$$
where $G_{ij}^l=\frac{\partial G_{i}^l}{\partial y^j}.$

\subsection{The mixed curvature of the Berwald connection:}
$$P_{ijk}^l=-G^{l}_{ijk}, \ \ \textrm{where}\ \ G_{ijk}^l=\frac{\partial G_{ij}^l}{\partial y^k}.$$

\subsection{An identity:}
\begin{equation}
\label{eq:4}
P^l_{ij}=-\frac{F}{2}l_m g^{kl}P_{ijk}^m
\end{equation}
\begin{Pf} Since
\begin{itemize}
\item $\displaystyle{Fl_m=\frac{\partial E}{\partial y^m},}$
\item $\displaystyle{\frac{\partial E}{\partial y^m}G^m=\frac{1}{2}y^k\frac{\partial E}{\partial x^k}}$, 
\item $\displaystyle{g_{mi}G^m=\frac{1}{2}\left(y^k\frac{\partial^2 E}{\partial y^i\partial x^k}-\frac{\partial E}{\partial x^i}\right)}$,
\item $\displaystyle{\frac{\partial} {\partial y^i}\left(\frac{\partial E}{\partial y^m}G^m\right)-g_{mi}G^m=\frac{\partial E}{\partial x^i}}$
\end{itemize}
we have
$$-Fl_mP^m_{ijk}=\frac{\partial E}{\partial y^m}G^m_{ijk}=$$
$$\frac{\partial}{\partial y^k}\left(\frac{\partial E}{\partial y^m}G^m_{ij}\right)-g_{mk}G_{ij}^m=\frac{\partial}{\partial y^k}\left(\frac{\partial}{\partial y^j} \left(\frac{\partial E}{\partial y^m}G^m_{i}\right)-g_{mj}G_i^m\right)-g_{mk}G_{ij}^m=$$
$$\frac{\partial}{\partial y^k}\left(\frac{\partial}{\partial y^j} \left(\frac{\partial} {\partial y^i}\left(\frac{\partial E}{\partial y^m}G^m\right)-g_{mi}G^m\right)-g_{mj}G^m_i\right)-g_{mk}G_{ij}^m=$$
$$\frac{\partial}{\partial y^k}\left(\frac{\partial}{\partial y^j} \left(\frac{\partial E}{\partial x^i}\right)-g_{mj}G^m_i\right)-g_{mk}G_{ij}^m=\frac{\partial}{\partial x^i}g_{jk}-2\C_{jmk}G_i^m-g_{mj}G_{ik}^m-g_{mk}G^m_{ij}=$$
$$2P_{ijk}=2g_{kl}P^l_{ij}\ \ \Rightarrow\ \ P^l_{ij}=-\frac{F}{2}l_m g^{kl}P_{ijk}^m$$
as was to be proved
\end{Pf}

\section{Conformality} 

\begin{Def}
Let $\tF$, $F\colon TM \to \mathbb{R}$ be Finsler metrics. They are conformally related if $\tF(x,y)=e^{\alpha(x)}F(x,y)$, where $\alpha \colon M\to \mathbb{R}$ is a function \emph{(}depending only on the position\emph{)}. 
\end{Def}

As an easy consequence of the conformality we have:
$$\tE=e^{2\alpha}E,\ \ \tilde{g}_{ij}=e^{2\alpha}g_{ij}\ \ \Rightarrow\ \ \tC_{ij}^k=\C_{ij}^k.$$
For the sake of simplicity let us use the following abbreviation: 
\begin{itemize}
\item $\displaystyle{\alpha_m=\frac{\partial \alpha}{\partial u^m}},\ \ m=1, \ldots, n.$
\end{itemize}
Using (\ref{eq:3})
$$\tG^l=G^l+y^m\alpha_my^l-Eg^{lm}\alpha_m.$$
According to its distinguished role let us introduce the gradient-type vector field
\begin{equation}
\label{eq:5}
X=X^l\frac{\partial}{\partial y^l},\ \ \textrm{where}\ \ X^l=g^{lm}\alpha_m\ \ \Rightarrow\ \ g_{lm}X^l=\alpha_m 
\end{equation}
and $X^l\alpha_l=X^lg_{lm}X^m=g(X,X)$ is just the Riemann-Finsler norm square of the vector field $X$.
Under this notation
\begin{equation}
\label{eq:6}\tG^l=G^l+y^m\alpha_m y^l-EX^l,
\end{equation}
\begin{equation}
\label{eq:7}
\tG_{i}^l=G_{i}^l+\alpha_iy^l+y^m\alpha_m\delta_i^l-\frac{\partial E}{\partial y^i}X^l-E\frac{\partial X^l}{\partial y^i}.
\end{equation}
According to its distinguished role let us introduce the (vector valued) one-form
\begin{equation}
\label{eq:8}
X^l_i=\frac{\partial X^l}{\partial y^i}.
\end{equation}
We have
\begin{equation}
\label{eq:9}
X_i^l=\frac{\partial X^l}{\partial y^i}=\frac{\partial g^{lm}}{\partial y^i}\alpha_m=-2\C_i^{lm}\alpha_m=\left\{\begin{array}{rl}
&-2g^{ms}\C_{is}^l\alpha_m=-2X^s\C_{is}^l\ \textrm{or, equivalently,}\\
&\\
&-2g^{ls}\C_{is}^m\alpha_m
\end{array}
\right.
\end{equation}
and the second line of formula (\ref{eq:9}) gives the {\bf{symmetry}} property
\begin{equation}
\label{eq:10}
g_{lk}X^l_i=-2\C_{ik}^m\alpha_m \ \ \Rightarrow\ \ g_{lk}X^l_i=g_{li}X^l_k
\end{equation}
because of the symmetry of the first Cartan tensor. Therefore we also have the following {\bf{cross-lifting}} formula
\begin{equation}
\label{eq:11}
g^{kl}X_k^mg_{mj}=X^l_j.
\end{equation}
Finally
\begin{equation}
\label{eq:12}
X^sX^l_s=-2X^sX^t\C_{st}^l.
\end{equation}
Using formula (\ref{eq:7}) 
\begin{equation}
\label{eq:13}
\tG_{ij}^l=G_{ij}^l+\alpha_i\delta_j^l+\alpha_j\delta_i^l-g_{ij}X^l-\frac{\partial E}{\partial y^i}X^l_j-\frac{\partial E}{\partial y^j}X^l_i-E\frac{\partial X_i^l}{\partial y^j}.
\end{equation}
According to its distinguished role let us introduce the quantity
$$X^l_{ij}=\frac{\partial X_i^l}{\partial y^j}=\frac{\partial^2 X^l}{\partial y^j \partial y^i}.$$ 

\begin{Prop} \emph{(Transformation formula for the Landsberg tensor)}
$$\tP_{ij}^l=P_{ij}^l-y^m\alpha_m\C_{ij}^l+EX_{i}^s\C_{js}^l+\frac{1}{2}\frac{\partial E}{\partial y^i}X^l_j+\frac{1}{2}\frac{\partial E}{\partial y^j}X^l_i+\frac{1}{2}g_{sj}X_i^sy^l+\frac{1}{2}EX_{ij}^l+\frac{1}{2}Eg^{lm}X^s_{im}g_{sj}.$$
\end{Prop}

The proof is a long straightforward calculation; see \cite{V3} and Hashiguchi \cite{H2}. In terms of the first Cartan tensor
$$g^{lm}X^s_{im}g_{sj}=g^{lm}\frac{\partial X_m^s}{\partial y^i}g_{sj}=\frac{\partial g^{lm}X_m^sg_{sj}}{\partial y^i}-\frac{\partial g^{lm}}{\partial y^i}X_m^sg_{sj}-g^{lm}X_m^s\frac{\partial g_{sj}}{\partial y^i}\stackrel{(\ref{eq:11})}{=}$$
$$X_{ij}^l+2\C^{lm}_i X_m^sg_{sj}-2g^{lm}X_m^sg_{sr}\C_{ji}^r\stackrel{(\ref{eq:10})}{=}$$
$$X_{ij}^l+2\C^{lm}_i X_j^sg_{sm}-2g^{lm}X_r^sg_{sm}\C_{ji}^r=X_{ij}^l+2\C^{l}_{is} X_j^s-2X_r^l\C_{ij}^r.$$
Therefore
\begin{equation}
\label{eq:14}
\tP_{ij}^l=P_{ij}^l-y^m\alpha_m\C_{ij}^l+\frac{1}{2}\frac{\partial E}{\partial y^i}X^l_j+\frac{1}{2}\frac{\partial E}{\partial y^j}X^l_i+
\frac{1}{2}g_{sj}X_i^sy^l+EX_{ij}^l+EX_{i}^s\C_{js}^l+E\C^{l}_{is} X_j^s-EX_r^l\C_{ij}^r.
\end{equation}

\begin{Rem} \emph{Some further computation shows that
$$X_{ij}^l+X_{i}^s\C_{js}^l+\C^{l}_{is} X_j^s-X_r^l\C_{ij}^r=-2X^s\nabla_s^v \C_{ij}^l$$
where $\nabla_s^v$ denotes the v-covariant derivative with respect to the Cartan connection. Therefore
$$\tP_{ij}^l=P_{ij}^l-y^m\alpha_m\C_{ij}^l+\frac{1}{2}\frac{\partial E}{\partial y^i}X^l_j+\frac{1}{2}\frac{\partial E}{\partial y^j}X^l_i+
\frac{1}{2}g_{sj}X_i^sy^l-2EX^s\nabla_s^v \C_{ij}^l$$
and the formula corresponds to formula (21) in \cite{V3} or formula (3.4C) in \cite{H2}; note that \cite{H2} involves an extra minus sign in the definition of $P_{ij}^l$ (formula (2.14*), page 35).} 
\end{Rem}

\section{Special conformal relationships I}

{\bf From now on we suppose that the Landsberg tensor satisfies the invariance property}
$$\tilde{P}_{ij}^l\alpha_l=P_{ij}^l\alpha_l$$
{\bf under the conformal change} $\tF=e^{\alpha}F$. As a direct consequence of formula (\ref{eq:14}) and the invariance property we have special expressions for $X_{ij}^l$ and the  contracted quantities $X_{ij}^l\alpha_l$ and $X^jX_{ij}^l$. 

\begin{Cor} If $\tP_{ij}^l=P_{ij}^l$ then
\begin{equation}
\label{eq:15}
EX_{ij}^l=y^m\alpha_m\C_{ij}^l-\frac{1}{2}\frac{\partial E}{\partial y^i}X^l_j-\frac{1}{2}\frac{\partial E}{\partial y^j}X^l_i-
\frac{1}{2}g_{sj}X_i^sy^l-EX_{i}^s\C_{js}^l-E\C^{l}_{is} X_j^s+EX_m^l\C_{ij}^m.
\end{equation}
\end{Cor}

\begin{Pf}
Equation (\ref{eq:15}) is a direct consequence of the transformation formula (\ref{eq:14}).
\end{Pf} 

\begin{Cor} If $\tP_{ij}^l\alpha_l=P_{ij}^l\alpha_l$ then
\begin{equation}
\label{eq:15a}
EX_{ij}^l\alpha_l=
\end{equation}
$$\left(y^m\alpha_m\C_{ij}^l-\frac{1}{2}\frac{\partial E}{\partial y^i}X^l_j-\frac{1}{2}\frac{\partial E}{\partial y^j}X^l_i-
\frac{1}{2}g_{sj}X_i^sy^l-EX_{i}^s\C_{js}^l-E\C^{l}_{is} X_j^s+EX_m^l\C_{ij}^m\right)\alpha_l$$
and
\begin{equation}
\label{eq:15b}
EX^jX_{ij}^l=
\end{equation}
$$X^j\left(y^m\alpha_m\C_{ij}^l-\frac{1}{2}\frac{\partial E}{\partial y^i}X^l_j-\frac{1}{2}\frac{\partial E}{\partial y^j}X^l_i-
\frac{1}{2}g_{sj}X_i^sy^l-EX_{i}^s\C_{js}^l-E\C^{l}_{is} X_j^s+EX_m^l\C_{ij}^m\right).$$
\end{Cor}

\begin{Pf}
Equation (\ref{eq:15a}) is a direct consequence of Corollary 1. Since the lowered second Cartan tensor $P_{ijk}$ is totally symmetric we have that 
$$P_{ij}^l\alpha_l=g^{lk}P_{ijk}\alpha_l \stackrel{(\ref{eq:5})}{=}X^kP_{ijk}=X^kP_{ikj}=X^kP^l_{ik}g_{lj}.$$
Using that $\tg_{ij}=e^{2\alpha}g_{ij}$, the same computation results in 
$$\tP_{ij}^l\alpha_l=\tilde{g}^{lk}\tP_{ijk}\alpha_l=\frac{1}{e^{2\alpha}}X^k\tP_{ijk}=\frac{1}{e^{2\alpha}}X^k\tP_{ikj}=\frac{1}{e^{2\alpha}}X^k\tP_{ik}^l\tilde{g}_{lj}=X^k\tP_{ik}^lg_{lj},$$
i.e. $\tP_{ij}^l\alpha_l=P_{ij}^l\alpha_l$ implies that
$$X^kP^l_{ik}g_{lj}=X^k\tP_{ik}^lg_{lj}\ \ \Rightarrow\ \ X^kP^l_{ik}=X^k\tP_{ik}^l.$$
This means, by Corollary 1, that formula (\ref{eq:15b}) holds.
\end{Pf}

\vspace{0.5cm}
\subsection{The first basic step} Equations (\ref{eq:15a}) and (\ref{eq:15b}) will be the key formulas to conclude the linear dependence of the vector fields 
$$X^l\frac{\partial}{\partial y^l}-\frac{X^l}{F}\frac{\partial F}{\partial y^l}C\ \ \ \textrm{and}\ \ \ X^sX_s^l\frac{\partial}{\partial y^l}\stackrel{(12)}{=}-2X^sX^t\C_{st}^l\frac{\partial}{\partial y^l}.$$
Since
$$X^l\frac{\partial F}{\partial y^l}=\frac{X^l}{F}\frac{\partial E}{\partial y^l}\stackrel{(5)}{=}\frac{g^{lm} \alpha_m}{F}\frac{\partial E}{\partial y^l}=\frac{y^m\alpha_m}{F}$$
we can also write that
$$X^l\frac{\partial}{\partial y^l}-\frac{X^l}{F}\frac{\partial F}{\partial y^l}C=X^l\frac{\partial}{\partial y^l}-\frac{y^m\alpha_m}{F^2}y^l\frac{\partial}{\partial y^l}=X^l\frac{\partial}{\partial y^l}-\frac{y^m\alpha_m}{2E}y^l\frac{\partial}{\partial y^l}.$$

\begin{Rem}
\emph{Note that the projected vector field 
$$X^l\frac{\partial}{\partial y^l}-\frac{X^l}{F}\frac{\partial F}{\partial y^l}C$$
is obviously tangential to the indicatrix hypersurface. The vector field $\displaystyle{X^sX_s^l\frac{\partial}{\partial y^l}}$ is also tangential to the indicatrix because of formula (\ref{eq:12}) and the basic properties of the first Cartan tensor}:
$$X^sX_s^l\frac{\partial F}{\partial y^l}=\frac{X^sX_s^l}{F}\frac{\partial E}{\partial y^l}\stackrel{(\ref{eq:9})}{=}\frac{X^s}{F}\left(-2X^t\mathcal{C}_{st}^l\frac{\partial E}{\partial y^l}\right),\ \ \textrm{where}\ \ \mathcal{C}_{st}^l\frac{\partial E}{\partial y^l}=g^{lm}\mathcal{C}_{mst}\frac{\partial E}{\partial y^l}=y^m\mathcal{C}_{mst}=0.$$
\end{Rem}

\vspace{0.5cm}
To conclude the linear dependency we use the substitution of $X^s$ systematically into the arguments of the difference tensor $\displaystyle{B_{ijk}^l=\tG_{ijk}^l-G_{ijk}^l}.$ Especially we prove the following lemma which is the generalization of Lemma 5 in \cite{V3} (page 22). 

\begin{Lem} If $\tP_{ij}^l\alpha_l=P_{ij}^l\alpha_l$ then
$$\frac{1}{3}X^kX^jB_{ijk}^l\alpha_l=\frac{1}{2}\left(X^l\alpha_lX^s_i\alpha_s-X^kX_k^j\alpha_j\alpha_i\right)+\frac{y^m\alpha_m}{4E}\left(X_j^tX^j\alpha_t\frac{\partial E}{\partial y^i}-(y^m\alpha_m)X_i^l\alpha_l\right)+$$
$$EX^jX^s_jX^pX^rQ_{pisr}$$
where
$$B_{ijk}^l=\tG_{ijk}^l-G_{ijk}^l$$
is the difference tensor of the mixed curveture of the Berwald connection and $Q_{pisr}=g_{pl}Q^l_{isr}$ is the lowered $vv$-curvature tensor of the Cartan connection.
\end{Lem}

\begin{Pf} Let us introduce the abbreviation
\begin{equation}
\label{eq:16}
B_{ij}^l=\tG_{ij}^l-G_{ij}^l\stackrel{(\ref{eq:13})}{=}\alpha_i\delta_j^l+\alpha_j\delta_i^l-g_{ij}X^l-\frac{\partial E}{\partial y^i}X^l_j-\frac{\partial E}{\partial y^j}X^l_i-EX_{ij}^l.
\end{equation}
Since
\begin{equation}
\label{eq:17}X^kX^jB_{ijk}^l\alpha_l=X^kX^j\frac{\partial B_{ij}^l}{\partial y^k}\alpha_l=X^k\frac{\partial}{\partial y^k}\left(X^jB_{ij}^l\alpha_l\right)-X^k X^j_k B_{ij}^l\alpha_l
\end{equation}
it is enough to compute the terms
$$X^k\frac{\partial}{\partial y^k}\left(X^jB_{ij}^l\alpha_l\right)\ \ \textrm{and}\ \ \ X^k X^j_k B_{ij}^l\alpha_l.$$
By some direct calculations
\begin{equation}
\label{eq:18}
X^jB_{ij}^l\alpha_l\stackrel{(\ref{eq:16})}{=}X^j\alpha_j\alpha_i-\frac{\partial E}{\partial y^i}X^jX^l_j\alpha_l-y^m\alpha_mX_i^l\alpha_l-EX^jX_{ij}^l\alpha_l
\end{equation}
because of
$$X^j\frac{\partial E}{\partial y^j}=g^{jm}\alpha_m\frac{\partial E}{\partial y^j}=y^m\alpha_m\ \ \ \textrm{and}\ \ \ X^jg_{ij}X^l=\alpha_iX^l.$$
Using formula (\ref{eq:15a}) 
$$EX^jX_{ij}^l\alpha_l=y^m\alpha_m\underline{X^j\C_{ij}^l\alpha_l}-\frac{1}{2}\frac{\partial E}{\partial y^i}X^jX^l_j\alpha_l-$$
$$\frac{1}{2}X^j\frac{\partial E}{\partial y^j}\underline{X^l_i\alpha_l}-\frac{1}{2}\underline{X^jg_{sj}X_i^s}y^l\alpha_l-EX^s_iX^j\C_{js}^l\alpha_l-E\C^{l}_{is}X^j X_j^s\alpha_l+EX_m^l\alpha_lX^j\C_{ij}^m\stackrel{(\ref{eq:5}), (\ref{eq:9})}{=}$$
$$-\frac{3}{2}y^m\alpha_mX_i^l\alpha_l-\frac{1}{2}\frac{\partial E}{\partial y^i}X^jX^l_j\alpha_l-E\C^{l}_{is}X^j X_j^s\alpha_l-EX^s_iX^j\C_{js}^l\alpha_l+EX_r^l\alpha_lX^j\C_{ij}^r.$$
On the other hand
$$X_r^l\alpha_lX^j\C_{ij}^r\stackrel{(\ref{eq:9})}{=}-\frac{1}{2}X_r^l\alpha_l X_i^r,\ \ X^s_iX^j\C_{js}^l\alpha_l\stackrel{(\ref{eq:9})}{=}-\frac{1}{2}X_i^sX_s^l\alpha_l,$$
$$\C^{l}_{is}X^j X_j^s\alpha_l\stackrel{(\ref{eq:10})}{=}-\frac{1}{2}g_{sr}X^r_i X^j X_j^s=-\frac{1}{2}g_{sj}X^r_i X^j X_r^s=-\frac{1}{2}\alpha_sX^r_i X_r^s$$
and, consequently,
\begin{equation}
\label{eq:19}
EX^jX_{ij}^l\alpha_l=-\frac{3}{2}y^m\alpha_mX_i^l\alpha_l-\frac{1}{2}\frac{\partial E}{\partial y^i}X^jX^l_j\alpha_l+\frac{1}{2}EX^l_i X_l^m\alpha_m.
\end{equation}
From equations (\ref{eq:18}) and (\ref{eq:19}) 
\begin{equation}
\label{eq:20}
X^jB_{ij}^l\alpha_l=X^j\alpha_j\alpha_i-\frac{1}{2}\frac{\partial E}{\partial y^i}X^jX^l_j\alpha_l+\frac{1}{2}y^l\alpha_l X_i^s\alpha_s-\frac{1}{2}EX^m_l X_i^l\alpha_m.
\end{equation}
Using the previous formula
\begin{equation}
\label{eq:full}
X^k\frac{\partial}{\partial y^k}\left(X^jB_{ij}^l\alpha_l\right)=X^k\left(X^j_k\alpha_j\alpha_i-\frac{1}{2}g_{ik}X^jX^l_j\alpha_l-\frac{1}{2}\frac{\partial E}{\partial y^i}X_k^jX^l_j\alpha_l-\frac{1}{2}\frac{\partial E}{\partial y^i}X^jX^l_{jk}\alpha_l\right)+
\end{equation}
$$X^k\left(\frac{1}{2}\alpha_k X_i^s\alpha_s+\frac{1}{2}y^l\alpha_l X_{ik}^s\alpha_s-\frac{1}{2}\frac{\partial E}{\partial y^k}X^m_l X_i^l\alpha_m-\frac{1}{2}EX^m_{lk} X_i^l\alpha_m-\frac{1}{2}EX^m_{l} X_{ik}^l\alpha_m\right)=$$
$$X^k X^j_k\alpha_j\alpha_i-\frac{1}{2}\alpha_iX^jX^l_j\alpha_l-\frac{1}{2}\frac{\partial E}{\partial y^i}X^kX_k^jX^l_j\alpha_l-\frac{1}{2}\frac{\partial E}{\partial y^i}X^jX^kX^l_{jk}\alpha_l+$$
$$\frac{1}{2}X^l\alpha_l X_i^s\alpha_s+\frac{1}{2}y^l\alpha_l X^k X_{ik}^s\alpha_s-\frac{1}{2}y^m\alpha_mX_l^s X^l_i\alpha_s-\frac{1}{2}EX^kX^m_{lk} X_i^l\alpha_m-\frac{1}{2}EX^m_{l} X^k X_{ik}^l\alpha_m=$$
$$\frac{1}{2}\alpha_iX^jX^l_j\alpha_l-\frac{1}{2}\frac{\partial E}{\partial y^i}X^kX_k^jX^l_j\alpha_l-\frac{1}{2}\frac{\partial E}{\partial y^i}X^j\underline{X^kX^l_{jk}\alpha_l}+$$
$$\frac{1}{2}X^l\alpha_l X_i^s\alpha_s+\frac{1}{2}y^l\alpha_l \underline{X^kX_{ik}^s\alpha_s}-\frac{1}{2}y^m\alpha_mX_l^s X^l_i\alpha_s-\frac{1}{2}EX_i^l\underline{X^kX^m_{lk} \alpha_m}-\frac{1}{2}\underline{E X^k X_{ik}^lX^m_l\alpha_m}.$$
Each term containing $X_{ij}^l$ means second order partial derivatives of $X^l$ with respect to $y$'s. To reduce the order of the partial differentiation in the first, the second and the third indicated terms we can directly use formula (\ref{eq:19}); for example (by replacing the free index $i$ with $j$ in (\ref{eq:19}))
\begin{equation}
\label{eq:21}
X^kX^l_{jk}\alpha_l\stackrel{(\ref{eq:19})}{=}\frac{1}{E}\left(-\frac{3}{2}y^m\alpha_m X_j^l\alpha_l-\frac{1}{2}\frac{\partial E}{\partial y^j}X^kX^l_k\alpha_l+\frac{1}{2}EX^l_j X_l^m\alpha_m\right)
\end{equation}
and we have similar expressions coming from the terms  
$$X^kX^m_{lk}\alpha_m=\frac{1}{E}\left(-\frac{3}{2}y^m\alpha_m X_l^t\alpha_t-\frac{1}{2}\frac{\partial E}{\partial y^l}X^kX^t_k\alpha_t+\frac{1}{2}EX^t_l X_t^m\alpha_m\right)$$
and
$$X^kX_{ik}^s\alpha_s=\frac{1}{E}\left(-\frac{3}{2}y^m\alpha_m X_i^l\alpha_l-\frac{1}{2}\frac{\partial E}{\partial y^i}X^kX^l_k\alpha_l+\frac{1}{2}EX^l_i X_l^m\alpha_m\right).$$
For the last indicated term 
$$EX^jX_{ij}^lX_l^t\alpha_t\stackrel{(\ref{eq:15b})}{=}$$
$$X^j\left(y^m\alpha_m\C_{ij}^l-\frac{1}{2}\frac{\partial E}{\partial y^i}X^l_j-\frac{1}{2}\frac{\partial E}{\partial y^j}X^l_i-
\frac{1}{2}g_{sj}X_i^sy^l-EX_{i}^s\C_{js}^l-E\C^{l}_{is} X_j^s+EX_m^l\C_{ij}^m\right)X_l^t\alpha_t\stackrel{(\ref{eq:9})}{=}$$
$$y^m\alpha_m\left(-\frac{1}{2}X_i^l\right)X_l^t\alpha_t-\frac{1}{2}\frac{\partial E}{\partial y^i}X^jX^l_jX_l^t\alpha_t-\frac{1}{2}y^m\alpha_m X^l_iX_l^t\alpha_t-
\frac{1}{2}\alpha_sX_i^sy^lX_l^t\alpha_t-$$
$$EX_{i}^s\left(-\frac{1}{2}X_s^l\right)X_l^t\alpha_t-EX^j\C^{l}_{is} X_j^sX_l^t\alpha_t+EX_m^l\left(-\frac{1}{2}X_i^m \right)X_l^t\alpha_t$$
because of 
$$\frac{1}{2}X^j\frac{\partial E}{\partial y^j}X^l_i\stackrel{(\ref{eq:5})}{=}\frac{1}{2}g^{jm}\alpha_m\frac{\partial E}{\partial y^j}X^l_i=\frac{1}{2}y^m\alpha_m X^l_i\ \ \textrm{and}\ \  y^lX_l^t=0;$$
see (\ref{eq:9}) and the basic property $y^l\C_{ijl}=0$ of the first Cartan tensor. Therefore
$$EX^jX_{ij}^lX_l^t\alpha_t=-y^m\alpha_m X_i^l X_l^t\alpha_t-\frac{1}{2}\frac{\partial E}{\partial y^i}X^jX^l_jX_l^t\alpha_t-EX^j\C^{l}_{is} X_j^sX_l^t\alpha_t$$
and we have by substituting the expressions of $X^kX^l_{jk}\alpha_l$ (see (\ref{eq:21})), $X^kX^m_{lk}\alpha_m$, $X^kX_{ik}^s\alpha_s$ and $EX^jX_{ij}^lX_l^t\alpha_t$ in 
formula (\ref{eq:full})
$$X^k\frac{\partial}{\partial y^k}\left(X^jB_{ij}^l\alpha_l\right)=\frac{1}{2}(X^kX_k^j\alpha_j)\alpha_i-\frac{1}{2}\frac{\partial E}{\partial y^i}X^jX^l_jX_l^t\alpha_t+
\frac{1}{2}X^l\alpha_lX_i^s\alpha_s-\frac{1}{2}y^m\alpha_m X_i^l X_l^t\alpha_t-$$
$$\frac{1}{2}\frac{\partial E}{\partial y^i}\left(-\frac{2}{E}y^m\alpha_mX^jX_j^l\alpha_l+\frac{1}{2}X^jX_j^lX_l^m\alpha_m\right)+$$
$$\frac{1}{2E}y^l\alpha_l\left( -\frac{3}{2}y^m\alpha_m X_i^l\alpha_l-\frac{1}{2}\frac{\partial E}{\partial y^i}X^t_jX^j \alpha_t+\frac{1}{2}EX^l_i X_l^m\alpha_m\right)-$$
$$\frac{1}{2}\left(-\frac{3}{2}y^m\alpha_m X_l^t\alpha_tX_i^l+\frac{1}{2}EX^l_i X_l^tX_t^m\alpha_m\right)-$$
$$\frac{1}{2}\left(-y^m\alpha_m X_i^l X_l^t\alpha_t-\frac{1}{2}\frac{\partial E}{\partial y^i}X^jX^l_jX_l^t\alpha_t-EX^j\C^{l}_{is} X_j^sX_l^t\alpha_t\right).$$
Now the formula has been free from the second order terms containing $X_{ij}^l$. 
On the other hand
$$X^kX_k^j\alpha_lB_{ij}^l\stackrel{(\ref{eq:16})}{=}X^kX_k^j\left(2\alpha_i\alpha_j-g_{ij}X^l\alpha_l-\frac{\partial E}{\partial y^i}X^l_j\alpha_l-\frac{\partial E}{\partial y^j}X^l_i\alpha_l-EX_{ij}^l\alpha_l\right).$$
To set the formula free from $X_{ij}^l$ we use formula (\ref{eq:15a}):
$$EX_{ij}^l\alpha_l=y^m\alpha_m\underline{\C_{ij}^l\alpha_l}-\frac{1}{2}\frac{\partial E}{\partial y^i}X^l_j\alpha_l-\frac{1}{2}\frac{\partial E}{\partial y^j}X^l_i\alpha_l-
\frac{1}{2}\underline{g_{sj}X_i^s}y^l\alpha_l-EX_{i}^s\C_{js}^l\alpha_l-E\C^{l}_{is} X_j^s\alpha_l+$$
$$EX_m^l\C_{ij}^m\alpha_l\stackrel{(\ref{eq:10})}=-\underline{g_{js}X^s_i}y^m\alpha_m-\frac{1}{2}\frac{\partial E}{\partial y^i}X^l_j\alpha_l-\frac{1}{2}\frac{\partial E}{\partial y^j}X^l_i\alpha_l+\frac{1}{2}EX_{i}^sg_{jm}X^m_s+\frac{1}{2}Eg_{sm}X^m_iX_j^s+$$
$$EX_m^l\C_{ij}^m\alpha_l$$
and, consequently,
$$X^kX_k^j\alpha_lB_{ij}^l=X^kX_k^j\left(2\alpha_i\alpha_j-g_{ij}X^l\alpha_l-\frac{1}{2}\frac{\partial E}{\partial y^i}X^l_j\alpha_l-\frac{1}{2}\frac{\partial E}{\partial y^j}X^l_i\alpha_l\right)+$$
$$X^kX_k^j\left(g_{js}X^s_iy^m\alpha_m-\frac{1}{2}EX_{i}^sg_{jm}X^m_s-\frac{1}{2}Eg_{sm}X^m_iX_j^s-EX_m^l\C_{ij}^m\alpha_l\right)=$$
$$2X^kX_k^j\alpha_i\alpha_j-g_{ij}X^kX_k^jX^l\alpha_l-\frac{1}{2}\frac{\partial E}{\partial y^i}X^kX_k^jX^l_j\alpha_l+g_{js}X^kX_k^j X^s_iy^m\alpha_m-$$
$$\frac{1}{2}EX_{i}^sg_{jm}X^kX_k^jX^m_s-\frac{1}{2}Eg_{sm}X^m_iX^kX_k^jX_j^s-EX_m^l\C_{ij}^mX^kX_k^j\alpha_l$$
because of
$$X^kX_k^j\frac{\partial E}{\partial y^j}=-2X^kX^t\C_{tk}^j\frac{\partial E}{\partial y^j}=-2X^kX^ty^m\C_{mtk}=0.$$
Using the symmetry property
$$EX_{i}^sg_{jm}X^kX_k^jX^m_s\stackrel{(\ref{eq:10})}{=}EX_{i}^sg_{sm}X^kX_k^jX^m_j$$
it follows that
$$X^kX_k^j\alpha_lB_{ij}^l=2X^kX_k^j\alpha_i\alpha_j-g_{ij}X^kX_k^jX^l\alpha_l-\frac{1}{2}\frac{\partial E}{\partial y^i}X^kX_k^jX^l_j\alpha_l+g_{js}X^kX_k^j X^s_iy^m\alpha_m-$$
$$EX_{i}^sg_{jm}X^kX_k^jX^m_s-EX_m^l\C_{ij}^mX^kX_k^j\alpha_l\stackrel{(\ref{eq:10})}{=}$$
$$2X^kX_k^j\alpha_i\alpha_j-g_{kj}X^kX_i^jX^l\alpha_l-\frac{1}{2}\frac{\partial E}{\partial y^i}X^kX_k^jX^l_j\alpha_l+g_{jk}X^kX_s^j X^s_iy^m\alpha_m-$$
$$EX_{i}^sg_{jk}X^kX_m^jX^m_s-EX_m^l\C_{ij}^mX^kX_k^j\alpha_l=$$
$$2X^kX_k^j\alpha_i\alpha_j-X_i^j\alpha_j X^l\alpha_l-\frac{1}{2}\frac{\partial E}{\partial y^i}X^kX_k^jX^l_j\alpha_l+X_s^j \alpha_j X^s_iy^m\alpha_m-$$
$$EX_{i}^sX_m^j\alpha_j X^m_s-EX_m^l\C_{ij}^mX^kX_k^j\alpha_l.$$
Finally
$$X^k\frac{\partial}{\partial y^k}\left(X^jB_{ij}^l\alpha_l\right)-X^kX_k^j\alpha_lB_{ij}^l=-\frac{3}{2}X^kX_k^j\alpha_j\alpha_i+\frac{3}{4E}y^m\alpha_mX_j^tX^j\alpha_t\frac{\partial E}{\partial y^i}+\frac{3}{2}X^l\alpha_lX^s_i\alpha_s-$$
$$\frac{3}{4E}(y^m\alpha_m)^2X_i^l\alpha_l+\frac{3}{2}EX^j\C^{l}_{is} X_j^sX_l^t\alpha_t+\frac{3}{4}EX_{i}^sX_m^j\alpha_j X^m_s,$$
where
$$\frac{3}{4}EX_{i}^sX_m^j\alpha_j X^m_s\stackrel{(\ref{eq:9})}{=}-\frac{3}{2}EX^t\C_{it}^sX_m^jX^m_s\alpha_j.$$
Therefore
$$X^k\frac{\partial}{\partial y^k}\left(X^jB_{ij}^l\alpha_l\right)-X^kX_k^j\alpha_lB_{ij}^l=-\frac{3}{2}X^kX_k^j\alpha_j\alpha_i+\frac{3}{4E}y^m\alpha_mX_j^tX^j\alpha_t\frac{\partial E}{\partial y^i}+\frac{3}{2}X^l\alpha_lX^s_i\alpha_s-$$
$$\frac{3}{4E}(y^m\alpha_m)^2X_i^l\alpha_l+\frac{3}{2}E\left(X^j\C^{l}_{is} X_j^sX_l^t\alpha_t-X^t\C_{it}^sX_m^jX^m_s\alpha_j\right),$$
where
$$X^j\C^{l}_{is} X_j^sX_l^t\alpha_t-X^t\C_{it}^sX_m^jX^m_s\alpha_j=2X^pX^jX_j^sX^rQ_{pisr}$$
as a straightforward computation shows\footnote{$$X^j\C^{l}_{is} X_j^sX_l^t\alpha_t-X^t\C_{it}^sX_m^jX^m_s\alpha_j\stackrel{(\ref{eq:5})}{=}X^j\C^{l}_{is} X_j^sX_l^t\alpha_t-X^t\C_{it}^sX_m^jX^m_sg_{jr}X^r\stackrel{(\ref{eq:10})}{=}$$
$$X^j\C^{l}_{is} X_j^sX_l^t\alpha_t-X^t\C_{it}^sX_r^jX^m_sg_{jm}X^r\stackrel{(\ref{eq:10})}{=}X^j\C^{l}_{is} X_j^sX_l^t\alpha_t-X^t\C_{it}^sX_r^jX^m_jg_{ms}X^r=$$
$$X^j\C^{l}_{is} X_j^sX_l^t\alpha_t-X^t\C_{mit}X_r^jX^m_jX^r=X^jX_j^s\left(\C^{l}_{is} X_l^t\alpha_t-X^t\C_{mit}X^m_s\right)\stackrel{(\ref{eq:5})}{=}X^jX_j^s\left(\C^{l}_{is} X_l^t\alpha_t-g^{tr}\alpha_r\C_{mit}X^m_s\right)=$$
$$X^jX_j^s\left(\C^{l}_{is} X_l^t\alpha_t-\alpha_r\C_{mi}^rX^m_s\right)=X^jX_j^s\left(\C^{l}_{is} X_l^t -\C_{mi}^tX^m_s\right)\alpha_t\stackrel{(\ref{eq:9})}{=}X^jX_j^s\left(\C^{l}_{is} \left(-2X^p\C_{pl}^t\right) -\C_{mi}^t\left(-2X^p\C_{ps}^m\right)\right)\alpha_t=$$
$$-2X^pX^jX_j^s\left(\C^{l}_{is} \C_{pl}^t -\C_{mi}^t\C_{ps}^m\right)\alpha_t=-2X^pX^jX_j^sQ^t_{ips}\alpha_t=-2X^pX^jX_j^sQ_{ips}^tX^rg_{rt}=-2X^pX^jX_j^sX^rQ_{ipsr}=$$
$$2X^pX^jX_j^sX^rQ_{pisr}\ \ \ \textrm{because of}\ \ \ Q_{ipsr}=-Q_{pisr}.$$}.
\end{Pf}

\begin{Cor}
$$\frac{1}{3}X^i_mX^mX^kX^jB_{ijk}^l\alpha_l=\frac{1}{2}\left(X^i_mX^mX^s_i\alpha_s\left(X^l\alpha_l-\frac{(y^m\alpha_m)^2}{2E}\right)-\left(X^kX_k^j\alpha_j\right)^2\right)+$$
$$EX^jX^s_jX^pX^rX^i_mX^mQ_{pisr},$$
where
$$X^i_mX^mX^s_i\alpha_s\left(X^l\alpha_l-\frac{(y^m\alpha_m)^2}{2E}\right)-\left(X^kX_k^j\alpha_j\right)^2=\textrm{the Gram determinant of the vector fields}$$
$$X^l\frac{\partial}{\partial y^l}-\frac{y^m\alpha_m}{2E}y^l\frac{\partial}{\partial y^l}\ \ \ \textrm{and}\ \ \ X^sX_s^l\frac{\partial}{\partial y^l}\ \ \textrm{with respect to the Riemann-Finsler metric} \ \ g_{ij}.$$
\end{Cor}
\section{Special conformal relationships II}
{\bf From now on we suppose that the mixed curvature of the Berwald connection satisfies the invariance property}
$$\tilde{P}_{ijk}^l\alpha_l=P_{ijk}^l\alpha_l$$
{\bf under the conformal change} $\tF=e^{\alpha}F$. Using subsection 1.6 this implies the invariance property
$$\tilde{P}_{ij}^l\alpha_l=P_{ij}^l\alpha_l$$
of the Landsberg tensor too. Since
$$\tilde{P}_{ijk}^l\alpha_l=P_{ijk}^l\alpha_l\ \ \Rightarrow\ \ B_{ijk}^l\alpha_l=0,$$
where $\displaystyle{B_{ijk}^l=\tG_{ijk}^l-G_{ijk}^l}$ is the difference tensor of the mixed curvatures, as a direct consequence of Corollary 3, we have 
$$0=\frac{1}{2E}\left(X^i_mX^mX^s_i\alpha_s\left(X^l\alpha_l-\frac{(y^m\alpha_m)^2}{2E}\right)-\left(X^kX_k^j\alpha_j\right)^2\right)+X^jX^s_jX^pX^rX^i_mX^mQ_{pisr}.$$
This means that if the tangent vector $v\in T_pM$ is of Finslerian length $1$ then 
$$F^2(v)=2E(v)=1\ \ \Rightarrow\ \ 0=1+\frac{X^jX^s_jX^pX^rX^i_mX^mQ_{pisr}}{\left(X^i_mX^mX^s_i\alpha_s\left(X^l\alpha_l-\frac{(y^m\alpha_m)^2}{2E}\right)-\left(X^kX_k^j\alpha_j\right)^2\right)}(v),$$
provided that we can divide by the Gram determinant. This is just the Gauss equation for the curvature of the indicatrix as a Riemannian submanifold\footnote{The Liouville vector field is the outer unit normal of the indicatrix with respect to $g_{ij}$. Since the components $\C_{ij}^k$ of the first Cartan tensor are the parameters of the L\'{e}vi-Civita connection and $y^j\C_{ij}^k=0$ it follows that the shape operator is working as the identity map and the normal curvature is constant $1$.}. To sum up: \emph{If the Gram determinant is not identically zero than we have zero sectional curvature of the indicatrix hypersurface} and the following theorem can be formulated.

\begin{Thm}  Let $M$ be a manifold of dimension $n\geq 3$ equipped with the Finslerian metric function $F$ and consider a function $\alpha \colon M \to \mathbb{R}$ satisfying the regularity condition $d_p\alpha\neq 0$ at a point $p\in M$.
If the sectional curvature of the indicatrix hypersurface at $p$ is positive and 
$$\tilde{P}_{ijk}^l\alpha_l=P^l_{ijk}\alpha_l$$
under the conformal change $\tilde{F}=e^{\alpha}F$, where $\alpha_l$ denotes the partial derivatives of $\alpha$ \emph{(}depending only on the position\emph{)}, then
the Gram determinant of the vector fields 
$$X^l\frac{\partial}{\partial y^l}-\frac{y^m\alpha_m}{2E}y^l\frac{\partial}{\partial y^l}\ \ \ \textrm{and}\ \ \ X^sX_s^l\frac{\partial}{\partial y^l}\ \ \textrm{with respect to the Riemann-Finsler metric} \ \ g_{ij}$$
must be zero and the vector fields are linearly dependent at all non-zero elements $v\in T_pM$.
\end{Thm}

Theorem 1 corresponds to Lemma 6 in \cite{V3} (page 28).

\subsection{The associated Riemannian metric} The Riemannian metric is derived from the conformal invariance property of the mixed curvature:
$$\tilde{P}_{ijk}^l\alpha_l=P_{ijk}^l\alpha_l\ \ \ \Rightarrow\ \ \ \frac{\partial^3 \tilde{G}^l}{\partial y^i \partial y^j \partial y^k}\alpha_l=\frac{\partial^3 G^l}{\partial y^i \partial y^j \partial y^k}\alpha_l$$
which means that the difference $G^l\alpha_l-\tilde{G}^l\alpha_l$ is quadratic in the tangent spaces. Therefore
\begin{equation}
\label{eq:22}
E_*:=EX^l\alpha_l\stackrel{(\ref{eq:6})}{=}\left(G^l-\tilde{G}^l\right)\alpha_l+(y^m\alpha_m)^2\ \ \textrm{is quadratic and}\ \ E_*\stackrel{(\ref{eq:5})}{=}EX^lg_{rl}X^r>0.
\end{equation}
Using the Riemannian energy function $E_*$ we will pay our attention to the associated Riemannian objects such as
$$g^*_{ij}=\frac{\partial E_*}{\partial y^i \partial y^j}\ \ \textrm{(Riemannian metric)},\ \ g_*^{ij}=(g^*_{ij})^{-1},\ \ X_*^l=g_*^{lk}\alpha_k\ \ \textrm{(Riemanian gradient)}.$$
Note that both $g_{ij}^*$ and its inverse depend only on the position, cf. Lemma 7 in  \cite{V3} (page 28).

\subsection{The second basic step: a Ricatti-type differential equation}

The following theorem is the key result to conclude a differential equation for the Finslerian energy $E$ along the lines
$$c(t)=v+tX_*(p),\ \ \textrm{where}\ \ v\in T_pM\ \ \ \textrm{and}\ \ \ v(\alpha)=0.$$
The proof will be presented in section 4.3. 
\begin{Thm}
The vector fields 
$$X^l\frac{\partial }{\partial y^l}\ \ (\textrm{the Riemann-Finsler gradient}),\ \ X_*^l\frac{\partial }{\partial y^l}\ \ (\textrm{the "vertically lifted" Riemannian gradient})$$
$$\textrm{and the Liouville vector field}\ \ C=y^l\frac{\partial }{\partial y^l}$$
form a linearly dependent system at all non-zero elements $v\in T_pM$, i.e. its Gram-determinant with respect to the metric $g_{ij}$ vanishes.
\end{Thm}

Theorem 2 corresponds to Lemma 9 in \cite{V3} (page 30). Using that the Gram-determinant of the vector fields vanishes we have:
\begin{equation}
0=\det \left(\begin{array}{ccc}
X^lX^kg_{lk}& X^lX_*^kg_{lk}&X^ly^kg_{lk}\\
&&\\
 X^lX_*^kg_{lk}&X^l_*X^k_*g_{lk}&X_*^ly^kg_{lk}\\
&&\\
X^ly^kg_{lk}&X_*^ly^kg_{lk}&y^ly^kg_{kl}
\end{array}\right)=\det \left(\begin{array}{ccc}
X^k\alpha_k& X_*^k\alpha_k&y^m\alpha_m\\
&&\\
 X_*^k\alpha_k&X^l_*X^k_* \displaystyle{\frac{\partial^2 E}{\partial y^k \partial y^l}}&X^k_* \displaystyle{\frac{\partial E}{\partial y^k}}\\
&&\\
y^m\alpha_m&X^k_* \displaystyle{\frac{\partial E}{\partial y^k}}&2E
\end{array}\right)=
\end{equation}
$$$$
$$X^k\alpha_k\left(2EX^l_*X^k_* \frac{\partial^2 E}{\partial y^k \partial y^l}-\left(X^k_* \frac{\partial E}{\partial y^k}\right)^2\right)-X^k_*\alpha_k\left(2EX_*^k\alpha_k-(y^m\alpha_m)X^k_* \frac{\partial E}{\partial y^k}\right)+$$
$$y^m\alpha_m\left((X_*^k\alpha_k)X^l_* \frac{\partial E}{\partial y^l}-(y^m\alpha_m)X^l_*X^k_* \frac{\partial^2 E}{\partial y^k \partial y^l}\right)\ \ \textrm{at all non-zero elements} \ \ v\in T_pM.$$
On the other hand
$$EX^k\alpha_k=E_*\ \ \textrm{(Riemannian energy)}\ \ \Rightarrow\ \ X^k\alpha_k=\frac{\ E_*}{E}$$
and the evaluation of the Gram-determinant along the line\footnote{It is the integral curve of the \emph{vertically lifted Riemannian gradient} vector field $\displaystyle X_*^l\frac{\partial}{\partial y^k}$ starting from $v$ in $T_pM$; note that the quantities $X_*^l$'s depend only on the position.} 
$$c(t)=v+tX_*(p),\ \ \textrm{where}\ \ v\in T_pM\ \ \textrm{and}\ \ v(\alpha)=0,$$
gives the differential equation
$$0=2E_*(c(t))(E\circ c)''(t)-E_*(c(t))\frac{\left((E\circ c)'(t)\right)^2}{(E\circ c)(t)}-2(X_*^k\alpha_k)^2(c(t))(E\circ c)(t)+$$
$$2(X_*^k\alpha_k)(c(t))(y^m\alpha_m)(c(t))(E\circ c)'(t)-(y^m\alpha_m)^2(c(t))(E\circ c)''(t).$$
Using the abbreviations
\begin{itemize}
\item $a^2=2E_*(v)$ the Riemannian norm square of the starting position $v\in T_pM$, where $v(\alpha)=0$,
\end{itemize}
\begin{itemize}
\item $b^2=2E_*(X_*(p))$ the Riemannian norm square of the Riemannian gradient vector field $X_*$ at the point $p$,
\end{itemize}
\begin{itemize}
\item $y(t)=E\circ c(t)$
\end{itemize}
we have that
\begin{itemize}
\item $2E_*(c(t))=a^2+t^2b^2$,
\end{itemize}
\begin{itemize}
\item $(X_*^k\alpha_k)(c(t))=(X_*^k\alpha_k)(p)=(g_*^{kl}\alpha_k\alpha_l)(p)=b^2,$
\end{itemize}
\begin{itemize}
\item $(y^m\alpha_m)(c(t))=v(\alpha)+t(X_*^m\alpha_m)(p)=tb^2\ \ \ \textrm{because of}\ \ v(\alpha)=0$
\end{itemize}
and the differential equation can be written into the following form 
\begin{equation}
\label{diffeqversion2d}
0=\left(a^2+t^2b^2(1-b^2)\right)y(t)y''(t)-\frac{1}{2}(a^2+t^2b^2)(y')^2(t)+2tb^4y(t)y'(t)-2b^4y^2(t).
\end{equation}
Using the substitution
$$z(t)=\frac{y'}{y}(t)$$
we obtain a Riccati-type differential equation
\begin{equation}
\label{diffeqversion2dricatti}
0=\left(a^2+t^2b^2(1-b^2)\right)z'(t)+\frac{1}{2}(a^2+t^2b^2(1-2b^2))z^2(t)+2tb^4z(t)-2b^4.
\end{equation}

\subsection{The proof of Theorem 2}

In order to prove Theorem 2 we need the relationship between the Riemannian metric $g^*_{ij}$ and the Riemann-Finsler metric $g_{ij}$. According to Theorem 1 we can use that the projected vector field 
\begin{equation}
\label{proj}
X^l\frac{\partial}{\partial y^l}-\frac{y^m\alpha_m}{2E}y^l\frac{\partial}{\partial y^l}
\end{equation}
and the vector field
\begin{equation}
\label{firstCartan}
X^sX_s^l\frac{\partial}{\partial y^l}
\end{equation}
are linearly dependent at all non-zero elements $v\in T_pM$. An easy calculation shows that the Riemann-Finsler norm square of the projected vector field (\ref{proj}) is
$$\left(X^l-\frac{y^m\alpha_m}{2E}y^l\right)\left(X^k-\frac{y^m\alpha_m}{2E}y^k\right)g_{kl}=\left(X^l-\frac{y^m\alpha_m}{2E}y^l\right)\left(\alpha_l-\frac{y^m\alpha_m}{2E}\frac{\partial E}{\partial y^l}\right)=$$
$$X^l\alpha_l-\frac{(y^m\alpha_m)^2}{2E}.$$
Therefore 
$$0\leq X^l\alpha_l-\frac{(y^m\alpha_m)^2}{2E}\ \ \Rightarrow\ \ 0\leq 2EX^l\alpha_l-(y^m\alpha_m)^2=2E_*-(y^m\alpha_m)^2\ \ \Rightarrow $$
\begin{equation}
\label{inequality}
\frac{(y^m\alpha_m)^2}{2E_*}\leq 1.
\end{equation}
The left hand side is homogeneous of degree $0$ and it attains its maximum at the (Riemannian) unit vector parallel to the (Riemannian) gradient $X_*(p)$. Since the Riemannian norm square of the Riemannian gradient vector field $X_*$ at the point $p$ is
$$(X_*^m\alpha_m)(p)=2E_*(X_*(p))\ \ \ \textrm{and}\ \ \ 2E_*=F_*^2$$
inequality (\ref{inequality}) shows that
$$b^2:=2E_*(X_*(p))\leq 1$$
and the equality 
$$0=2E_*(v)-(y^m\alpha_m)^2(v)$$
occours (i.e. the projected vector field vanishes at some nonzero element $v\in T_pM$) if and only if 
$$\frac{y^m\alpha_m}{F_*}(X_*(p))=\sup_{v\neq {\bf{0}}}\frac{y^m\alpha_m}{F_*}(v)=1.$$

\begin{Lem} \emph{(The zeros of the projected vector field)} The projected vector field \emph{(\ref{proj})}
\begin{itemize}
\item  does not vanish at any $v\neq \pm X_*(p)$ in $T_pM$,
\item vanishes at $\pm X_*(p)$ if and only if $\ \displaystyle{b^2:=2E_*(X_*(p))=1},$
i.e. the Riemannian gradient $X_*$ is of unit Riemannian lenght at the point $p$. 
\end{itemize}
\end{Lem}

The result corresponds to Lemma 8 in \cite{V3} (page 30). Using Theorem 1 and Lemma 2 we can write that
\begin{equation}
\label{eq:29}
X^l_sX^s=\theta\left(X^l-\frac{y^m\alpha_m}{2E}y^l\right)\ \ \textrm{for all}\ \ v\in T_pM\setminus\{\pm X_*(p)\}.
\end{equation}
According to the distinguished role of the projected vector field let us introduce the abbreviations
\begin{itemize} 
\item $\displaystyle{\eta^l=X^l-\frac{y^m\alpha_m}{2E}y^l\ \ \Rightarrow\ \ \eta:=\eta^l\alpha_l=X^l\alpha_l-\frac{(y^m\alpha_m)^2}{2E}.}$
\end{itemize}
From (\ref{eq:29})
\begin{equation}
\label{eq:30}
X^l_sX^s\alpha_l=\theta \left(X^l\alpha_l-\frac{(y^m\alpha_m)^2}{2E}\right).
\end{equation}
Differentiating the left hand side by $y$'s
$$\frac{\partial X^l_sX^s}{\partial y^j}\alpha_l=X^l_{js}X^s\alpha_l+X^l_sX_j^s\alpha_l\stackrel{(\ref{eq:21})}{=}-\frac{1}{2E}\left(3y^m\alpha_m X_j^l\alpha_l+\frac{\partial E}{\partial y^j}X^kX^l_k\alpha_l-EX^l_j X_l^m\alpha_m\right)+X^l_sX_j^s\alpha_l=$$
$$-\frac{1}{2E}\left(3y^m\alpha_m X_j^l\alpha_l+\frac{\partial E}{\partial y^j}X^kX^l_k\alpha_l\right)+\frac{3}{2}X^l_sX_j^s\alpha_l.$$
By the symmetry property
\begin{equation}
\label{eq:31}
X_j^l\alpha_l=X_j^lg_{lr}X^r\stackrel{(\ref{eq:10})}=X_r^lg_{lj}X^r\stackrel{(\ref{eq:29})}{=}\theta\left(\alpha_j-\frac{y^m\alpha_m}{2E}\frac{\partial E}{\partial y^j}\right).
\end{equation}
In a similar way,
$$X^kX^l_k\alpha_l=X^k\theta\left(\alpha_k-\frac{y^m\alpha_m}{2E}\frac{\partial E}{\partial y^k}\right)=\theta\left(X^k\alpha_k-\frac{(y^m\alpha_m)^2}{2E}\right)=\theta \eta$$
because of
$$X^k\frac{\partial E}{\partial y^k}=g^{kr}\alpha_r\frac{\partial E}{\partial y^k}=y^r\alpha_r$$
and 
$$X^l_sX_j^s\alpha_l\stackrel{(\ref{eq:31})}{=}\theta\left(\alpha_s-\frac{y^m\alpha_m}{2E}\frac{\partial E}{\partial y^s}\right)X_j^s=\theta \alpha_sX_j^s$$
because of
$$\frac{\partial E}{\partial y^s} X_j^s\stackrel{(\ref{eq:9})}{=}-2X^t\C_{tj}^s\frac{\partial E}{\partial y^s}=-2X^tg^{sr}\C_{tjr}\frac{\partial E}{\partial y^s}=-2X^ty^r\C_{tjr}=0.$$
Therefore
\begin{equation}
\label{eq:32}
X^l_sX_j^s\alpha_l=\theta \alpha_sX_j^s\stackrel{(\ref{eq:31})}{=}\theta^2\left(\alpha_j-\frac{y^m\alpha_m}{2E}\frac{\partial E}{\partial y^j}\right).
\end{equation}
Using the previous formulas
\begin{equation}
\label{eq:33}
\frac{\partial X^l_sX^s}{\partial y^j}\alpha_l=-\frac{1}{2E}\left(3y^m\alpha_m \theta \alpha_j+\left(\theta \eta-3\theta\frac{(y^m\alpha_m)^2}{2E}\right)\frac{\partial E}{\partial y^j}\right)+
\end{equation}
$$\frac{3}{2}\theta^2\left(\alpha_j-\frac{y^m\alpha_m}{2E}\frac{\partial E}{\partial y^j}\right).$$
Now we are going to differentiate the right hand side of (\ref{eq:30}):
\begin{equation}
\label{eq:34}
\frac{\partial \theta}{\partial y^j}\left(X^l\alpha_l-\frac{(y^m\alpha_m)^2}{2E}\right)+\theta \left(X^l_j\alpha_l-\frac{4E(y^m\alpha_m)\alpha_j-2(y^m\alpha_m)^2\frac{\partial E}{\partial y^j}}{(2E)^2}\right)=
\end{equation}
$$\frac{\partial \theta}{\partial y^j}\left(X^l\alpha_l-\frac{(y^m\alpha_m)^2}{2E}\right)+\theta \left(\theta\left(\alpha_j-\frac{y^m\alpha_m}{2E}\frac{\partial E}{\partial y^j}\right)-\frac{4E(y^m\alpha_m)\alpha_j-2(y^m\alpha_m)^2\frac{\partial E}{\partial y^j}}{(2E)^2}\right)=$$
$$\frac{\partial \theta}{\partial y^j}\eta+\left(\theta^2-\theta\frac{y^m\alpha_m}{E}\right)\alpha_j-\left(\theta^2\frac{y^m\alpha_m}{2E}-2\theta \frac{(y^m\alpha_m)^2}{(2E)^2}\right)\frac{\partial E}{\partial y^j}.$$
Comparing the partial derivatives (\ref{eq:33}) and (\ref{eq:34}) it follows that  
\begin{equation}
\label{eq:35}
\frac{\partial \theta}{\partial y^j}\eta=\left(-\theta\frac{\eta}{2E}+\theta \frac{(y^m\alpha_m)^2}{(2E)^2}-\frac{1}{2}\theta^2\frac{y^m\alpha_m}{2E}\right)\frac{\partial E}{\partial y^j}+\left(\frac{1}{2}\theta^2-\frac{1}{2}\theta\frac{y^m\alpha_m}{E}\right)\alpha_j.
\end{equation}
Now we are in the position to compute a detailed relationship between the Riemannian metric $g^*_{ij}$ and the Riemann-Finsler metric $g_{ij}$
as follows: $E_*=EX^l\alpha_l\ \ \Rightarrow$
$$g^*_{ij}=\frac{\partial^2 E_*}{\partial y^i \partial y^j}=\frac{\partial}{\partial y^j}\left(\frac{\partial E}{\partial y^i}X^l\alpha_l+EX_i^l\alpha_l\right)=\frac{\partial}{\partial y^j}\left(\frac{\partial E}{\partial y^i}X^l\alpha_l+E\theta\left(\alpha_i-\frac{y^m\alpha_m}{2E}\frac{\partial E}{\partial y^i}\right)\right)\stackrel{(\ref{eq:31})}{=}$$
$$g_{ij}X^l\alpha_l+\frac{\partial E}{\partial y^i}\theta\left(\alpha_j-\frac{y^m\alpha_m}{2E}\frac{\partial E}{\partial y^j}\right)+\frac{\partial E}{\partial y^j}\theta\left(\alpha_i-\frac{y^m\alpha_m}{2E}\frac{\partial E}{\partial y^i}\right)+$$
$$E\frac{\partial \theta}{\partial y^j}\left(\alpha_i-\frac{y^m\alpha_m}{2E}\frac{\partial E}{\partial y^i}\right)-E\theta\frac{y^m\alpha_m}{2E}g_{ij}-$$
$$E\theta\frac{2E\alpha_j-2(y^m\alpha_m)\frac{\partial E}{\partial y^j}}{(2E)^2}\frac{\partial E}{\partial y^i}.$$
Using (\ref{eq:35}) a straightforward calculation shows that
\begin{equation}
\label{eq:36}
g^*_{ij}=\left(X^l\alpha_l-\theta\frac{y^m\alpha_m}{2}\right)g_{ij}+\theta\frac{y^m\alpha_m}{4E}\left(\frac{y^m\alpha_m}{2\eta}\left(\theta-\frac{y^m\alpha_m}{E}\right)-1\right)\frac{\partial E}{\partial y^i}\frac{\partial E}{\partial y^j}+
\end{equation}
$$\theta\frac{E}{2\eta}\left(\theta-\frac{y^m\alpha_m}{E}\right)\alpha_i\alpha_j+\frac{\theta}{2}\left(1+\frac{(y^m\alpha_m)^2}{2E\eta}-\theta\frac{y^m\alpha_m}{2\eta}\right)\left(\alpha_i\frac{\partial E}{\partial y^j}+\alpha_j\frac{\partial E}{\partial y^i}\right).$$
Therefore
\begin{equation}
\label{eq:37}
g_{ij}^*\in \mathcal{L}\left(g_{ij}, \frac{\partial E}{\partial y^i}\frac{\partial E}{\partial y^j}, \alpha_i\alpha_j, \alpha_i\frac{\partial E}{\partial y^j}+\alpha_j\frac{\partial E}{\partial y^i}\right),
\end{equation}
where the linear combination contains functions as coefficients of the symmetric terms
$$g_{ij}, \ \ \frac{\partial E}{\partial y^i}\frac{\partial E}{\partial y^j}, \ \ \alpha_i\alpha_j,\ \ \alpha_i\frac{\partial E}{\partial y^j}+\alpha_j\frac{\partial E}{\partial y^i}.$$ 
Formula (\ref{eq:37}) implies the linear dependence of the Riemannian gradient $X_*$, the Riemann-Finsler gradient $X$ and the Liouville vector field $C$ because of
$$\alpha_j=X_*^ig_{ij}^*\in \mathcal{L}\left(X_*^ig_{ij}, X_*^i\frac{\partial E}{\partial y^i}\frac{\partial E}{\partial y^j}, X_*^i \alpha_i\alpha_j, X_*^i\left(\alpha_i\frac{\partial E}{\partial y^j}+\alpha_j\frac{\partial E}{\partial y^i}\right)\right)\ \ \Rightarrow\ \ $$
$$X^l=g^{lj}\alpha_j\in \mathcal{L}\left(X_*^l, X_*^i\frac{\partial E}{\partial y^i}y^l, X_*^i \alpha_i X^l, X_*^i\alpha_i y^l+X_*^i\frac{\partial E}{\partial y^i}X^l\right)$$
as was to be proved.

\subsection{What about the constant $b^2$?}

The solution of differential equation (\ref{diffeqversion2dricatti}) seems to be hard in general. In what follows we present an essential simplification in case of dimension $n\geq 3$. The result corresponds to Lemma 10 in \cite{V3} (page 32). Keeping the singular solutions in mind note that the forthcoming argument can be also used under some $y$ - locality:
\begin{itemize}
\item [(RP)] the Finslerian energy function $E$ is of class $\mathcal{C}^2$ at $y=X_*(p)$. 
\end{itemize}
\begin{Thm} $\displaystyle{b^2:=2E_*(X_*(p))=1}$.
\end{Thm}

\begin{Pf}
Suppose now that $b^2 < 1$. Lemma 2 shows that the projected vector field (\ref{proj}) has no zeros, i.e. its Finslerian norm square is strictly positive:
$$\eta=X^l\alpha_l-\frac{(y^m\alpha_m)^2}{2E}>0.$$
By (\ref{eq:29})
$$\theta=\frac{X^l_sX^s\alpha_l}{\eta}=\frac{1}{\eta}\frac{\partial X^l\alpha_l}{\partial y^s}X^s.$$
This means\footnote{Note that (RP) does not imply automatically that $X_s^l$ and the further $y$ - derivatives exist at $X_*(p)$ but we can introduce the contracted terms $X_s^l\alpha_l$ or $X_{js}^l \alpha_l$ up to order $2$ by the formulas
$$X_s^l\alpha_l=\frac{\partial X^l\alpha_l}{\partial y^s}\ \ \textrm{and}\ \ X_{js}^l\alpha_l=\frac{\partial^2 X^l\alpha_l}{\partial y^j\partial y^s}$$
because of $X^l\alpha_l=E_*/E$ and property (RP).} that $\theta$ is of class $\mathcal{C}^1$ at $X_*(p)$ because of $X^l\alpha_l=E_*/E$ and property (RP). 
In what follows we will use equation (\ref{eq:36}) at $v:=X_*(p)$. 
For the sake of simplicity let us introduce the following abbreviations
\begin{itemize}
\item $\displaystyle{A:=\left(X^l\alpha_l-\theta\frac{y^m\alpha_m}{2}\right)}$,
\end{itemize}
\begin{itemize}
\item $\displaystyle{P:=\theta\frac{y^m\alpha_m}{4E}\left(\frac{y^m\alpha_m}{2\eta}\left(\theta-\frac{y^m\alpha_m}{E}\right)-1\right)}$,
\end{itemize}
\begin{itemize}
\item $\displaystyle{R:=\theta\frac{E}{2\eta}\left(\theta-\frac{y^m\alpha_m}{E}\right)},$
\end{itemize}
\begin{itemize}
\item $\displaystyle{Q:=\frac{\theta}{2}\left(1+\frac{(y^m\alpha_m)^2}{2E\eta}-\theta\frac{y^m\alpha_m}{2\eta}\right)}$
\end{itemize}
for the coefficients in formula (\ref{eq:36}):
$$g^*_{ij}=Ag_{ij}+P\frac{\partial E}{\partial y^i}\frac{\partial E}{\partial y^j}+R\alpha_i\alpha_j+Q\left(\alpha_i\frac{\partial E}{\partial y^j}+\alpha_j\frac{\partial E}{\partial y^i}\right);$$
see the notations in our main reference work \cite{V3}, page 32, formula (43). An easy computation shows that
\begin{equation}
\label{eq:38}
2EP+y^m\alpha_m Q=0\ \ \textrm{and}\ \ \ 2EQ+Ry^m\alpha_m=E\theta.
\end{equation}
We have
$$X^k=g^{kj}\alpha_j=g^{kj}X_*^ig^*_{ij}\stackrel{(\ref{eq:36})}{=}AX_*^k+\left(PX_*^i\frac{\partial E}{\partial y^i}+Q X_*^i\alpha_i\right)y^k+$$
$$\left(QX_*^i\frac{\partial E}{\partial y^i}+RX_*^i\alpha_i\right)X^k.$$
Evaluating both side at $v:=X_*(p)$:
\begin{equation}
\label{CX}
\left(1-E\theta\right)(v)X^k(v)=A(v)v^k\ \ \ (k=1, \ldots, n)
\end{equation}
because of (\ref{eq:38}) and the homogenity property
$$v^i\frac{\partial E}{\partial y^i}(v)=2E(v).$$
Observe that $A(v)=0$ contradicts to equation (\ref{eq:36}). Indeed, {\bf since the dimension is at least $3$} we can choose a tangent vector $w\in T_pM$ to eliminate each term of 
$$\frac{\partial E}{\partial y^i}\frac{\partial E}{\partial y^j}, \ \ \alpha_i\alpha_j\ \ \textrm{and}\ \ \alpha_i\frac{\partial E}{\partial y^j}+\alpha_j\frac{\partial E}{\partial y^i}$$
at $v\in T_pM$. Therefore
$$w^iw^jg^*_{ij}(p)=A(v)w^iw^jg_{ij}(v).$$
This means that the main coefficient $A(v)$ must be positive because the Riemannian metric $g^*_{ij}$ (depending only on the position) is positive definite. Since $A(v)\neq0$ we have by (\ref{CX}) that $\left(1-E\theta\right)(v)\neq 0$ and
$$X^k(v):=\frac{A(v)}{\left(1-E\theta\right)(v)}v^k \ \ (k=1, \ldots, n)\ \ \Rightarrow \ \ X(v)\ || \ C(v),$$
i.e. the projected vector field (\ref{proj}) vanishes at $v:=X_*(p)$ which contradicts to $b^2 < 1$ in the sense of Lemma 2. 
\end{Pf}

\subsection{The initial condition}

\begin{Thm} $\displaystyle{\frac{y'(0)}{y(0)}=\frac{K(p)}{F_*(v)}}$
for some real constant $K(p)$.

\end{Thm}

In what follows {\bf we use again that the dimension is at least $3$}. This means that the linear subspace 
\begin{equation}
\label{subspace}
T_pN:=\{v\in T_pM\ | \ v(\alpha)=0 \}
\end{equation}
is of dimension at least $2$ and $T_pN\setminus \{{\bf 0}\}$ is connected.

\begin{Thm} Both  $\displaystyle{X^l\alpha_l(v)}$ and $\displaystyle{(FX_s^lX^s \alpha_l)(v)}$
are independent of the choice $v\in T_pN$.
\end{Thm}

\begin{Pf} Let $N$ be the level hypersurface of the function $\alpha$ such that $p\in N$ and consider a tangent vector $v$ at the point $p$. Formula (\ref{eq:36}) shows that 
$$g^*_{ij}(v)=X^l(v)\alpha_lg_{ij}(v)$$
and Knebelman's theorem for conformally related Riemann-Finsler metrics says that the scale function $X^l\alpha_l$ restricted to $T_pN$ is constant. We can get the same result by some direct calculations too: let $w\in T_pM$ be an arbitrary vector such that $w(\alpha)=0$ (geometrically $w$ is tangential to the level hypersurface of $\alpha$ passing through the point $p$). Then
$$w^i\frac{\partial X^l \alpha_l}{\partial y^i}(v)=w^iX_i^l(v) \alpha_l\stackrel{(\ref{eq:31})}{=}\theta(v)\left(w^i\alpha_i-\frac{v^m\alpha_m}{2E(v)}w^i\frac{\partial E}{\partial y^i}(v)\right)=0,$$
i.e. $\displaystyle{X^l\alpha_l}$ is constant along the subspace (\ref{subspace}). On the other hand
$$w^i\frac{\partial X_s^l X^s \alpha_l}{\partial y^i}(v)=w^iX_{si}^{l}(v)X^s(v)\alpha_l+w^iX_{s}^{l}(v)X_i^s(v)\alpha_l,$$
where, by formula (\ref{eq:32}), the second term vanishes:
\begin{equation}
\label{incon}
w^iX_{s}^{l}(v)X_i^s(v)\alpha_l=\theta^2(v)\left(w^i\alpha_i-\frac{v^m\alpha_m}{2E(v)}w^i\frac{\partial E}{\partial y^i}(v)\right)=0.
\end{equation}
Formula (\ref{eq:19}) and (\ref{incon}) show that
$$w^iX_{si}^{l}(v)X^s(v)\alpha_l=-\frac{1}{2E(v)}w^i\frac{\partial E}{\partial y^i}(v)X_{s}^{l}(v)X^s(v)\alpha_l.$$
Therefore
$$w^i\frac{\partial X_s^l X^s \alpha_l}{\partial y^i}(v)+\frac{1}{2E(v)}w^i\frac{\partial E}{\partial y^i}(v)X_{s}^{l}(v)X^s(v)\alpha_l=0\ \ \Rightarrow$$
$$F(v)w^i\frac{\partial X_s^l X^s \alpha_l}{\partial y^i}(v)+w^i\frac{\partial F}{\partial y^i}(v)X_{s}^{l}(v)X^s(v)\alpha_l=0$$
because of $F^2=2E$. Finally
$$w^i\frac{\partial F X_s^l X^s \alpha_l}{\partial y^i}(v)=0$$
as was to be proved. 
\end{Pf}

\vspace{0.5cm}
Theorem 5 corresponds to Lemma 11 in  \cite{V3} (page 33). Since
$$E=\frac{E_*}{X_l\alpha_l}$$
we have that
$$y'(0)=X_*^i(p)\frac{\partial E}{\partial y^i}(v)=X_*^i(p)\frac{\frac{\partial E_*}{\partial y^i}(v) (X^l\alpha_l)(v)-E_*(v)X_i^l(v)\alpha_l}{(X_l\alpha_l)^2(v)},$$
where
$$X_*^i(p)\frac{\partial E_*}{\partial y^i}(v)=g_*^{ij}(p)\alpha_j\frac{\partial E_*}{\partial y^i}(v)=(y^j\alpha_j)(v)=v^j\alpha_j=0$$
because $X_*$ is the Riemannian gradient of $\alpha$ and $v(\alpha)=0$. Therefore
$$y'(0)=-X_*^i(p)\frac{E_*(v)X_i^l(v)\alpha_l}{(X_l\alpha_l)^2(v)}\stackrel{(\ref{eq:31})}{=}-E_*(v)\theta(v)\frac{X_*^i(p)\left(\alpha_i-\frac{v^m\alpha_m}{2E(v)}\frac{\partial E}{\partial y^i}(v)\right)}{(X_l\alpha_l)^2(v)}=$$
$$-E_*(v)X_*^i(p)\alpha_i \frac{\theta(v)}{(X_l\alpha_l)^2(v)}\stackrel{(\ref{eq:30})}{=}-E_*(v)X_*^i(p)\alpha_i \frac{X_{s}^l(v)X^s(v)\alpha_l}{(X_l\alpha_l)^3(v)}$$
because of $v(\alpha)=0$. Finally
$$y'(0)=-E_*(v)X_*^i(p)\alpha_i \frac{F(v)X_{s}^l(v)X^s(v)\alpha_l}{F(v)(X_l\alpha_l)^3(v)}\stackrel{(\ref{eq:22})}{=}$$
$$-E_*(v)X_*^i(p)\alpha_i \frac{F(v)X_{s}^l(v)X^s(v)\alpha_l}{\underbrace{F(v)(X_l\alpha_l)^{1/2}(v)}_{F_*(v)} (X_l\alpha_l)^{5/2}(v)}=$$
$$-\frac{1}{2}F_*(v)X_*^i(p)\alpha_i \frac{F(v)X_{s}^l(v)X^s(v)\alpha_l}{(X_l\alpha_l)^{5/2}(v)}=-\frac{1}{2}F_*(v)\frac{F(v)X_{s}^l(v)X^s(v)\alpha_l}{(X_l\alpha_l)^{5/2}(v)}$$
because of
$$b^2=2E_*(X_*(p))=X_*^i(p)\alpha_i=1$$
(see section 4.4). On the other hand
$${y(0)}=E(v)=\frac{E_*(v)}{(X_l\alpha_l)(v)}$$
and, by Theorem 5, we can write that
\begin{equation}
\label{initial}
\frac{y'(0)}{y(0)}=\frac{K(p)}{F_*(v)}
\end{equation}
for some real constant $K(p)$.

\subsection{The solution of the differential equation}

Using section 4.4 (Theorem 3) we can reduce differential equation (\ref{diffeqversion2dricatti}) to the following simple form
\begin{equation}
\label{diffeqversion3dricatti}
0=a^2z'(t)+\frac{1}{2}(a^2-t^2)z^2(t)+2tz(t)-2.
\end{equation}
It can be directly seen that the function
\begin{equation}
\label{solution}
z(t):=2\frac{2t+K(p)F_*(v)}{2t^2+tK(p)F_*(v)+4E_*(v)}
\end{equation}
is the solution of the problem satisfying
$$z(0)=\frac{K(p)}{F_*(v)}.$$
Since the solution should be defined for all real parameters we have that
\begin{equation}
\label{asanovcondition1}
-4 < K(p) < 4.
\end{equation}
Integrating (\ref{solution}) 
\begin{equation}
\label{main}
E(v+tX_*(p))=K^*(p)\left(F_*^2(v)+K(p)F_*(v)\frac{t}{2}+t^2\right)e^{2A(v,t)},
\end{equation}
where
$$A(v,t)=\frac{\ \ K(p)}{\sqrt{16-K^2(p)}}\left(\arctan \frac{1}{\sqrt{16-K^2(p)}}\left(\frac{4t}{F_*(v)}+K(p)\right)-\arctan \frac{K(p)}{\sqrt{16-K^2(p)}}\right).
$$

\section{The comparison of the metrics}

Following Asanov \cite{A2} (section 2) we show that the metric (\ref{main}) belongs to the class of Finsleroid-Finsler metrics up to conformality. 

\begin{Thm}
The metric given by formula \emph{(\ref{main})} is conformal to a Finsleroid-Finsler metric.
\end{Thm}

\begin{Pf}
The general form of Finsleroid-Finsler metrics is given by
\begin{equation}
\label{asanovmetric1}
F=e^{G\Phi/2}\sqrt{b^2+gqb+q^2},
\end{equation}
where
\begin{itemize}
\item $b=b_i({\bf x})y^i$ (the Finsleroid axis $1$ - form),
\end{itemize}
\begin{itemize}
\item $a_{ij}({\bf x})y^iy^j$ (the Riemannian metric) and $a^{ij}b_ib_j=1$,
\end{itemize}
\begin{itemize}
\item $q=\sqrt{r_{ij}({\bf x})y^iy^j}$, where $r_{ij}=a_{ij}-b_ib_j$, 
\end{itemize}
\begin{itemize}
\item $g=g({\bf x})$ and $-2 < g({\bf x})< 2$ (the Finsleroid charge),
\end{itemize}
\begin{itemize}
\item $\displaystyle{h=\sqrt{1-\frac{\ g^2}{4}}},$
\end{itemize}
\begin{itemize}
\item $G=g/h,$
\end{itemize}
\begin{itemize} 
\item $\displaystyle{\Phi=\left\{\begin{array}{rl}
&+\frac{\pi}{2}+\arctan \frac{G}{2}-\arctan \frac{q+\frac{g}{2}b}{hb}\ \ \textrm{if} \ \ b > 0\\
&\\
&-\frac{\pi}{2}+\arctan \frac{G}{2}-\arctan \frac{q+\frac{g}{2}b}{hb}\ \ \textrm{if} \ \ b < 0.
\end{array}
\right.}$
\end{itemize}
The common limit of the right hand sides as $b\to 0$ is $\displaystyle{\arctan \frac{G}{2}}$. Let us introduce the function
$$f(b):=\arctan \frac{q+\frac{g}{2}b}{hb};$$
then
\begin{equation}
\label{f}
\lim_{b\to \pm \infty}f(b)=\arctan \frac{g}{2h}=\arctan \frac{G}{2}\ \ \textrm{and}\ \ f'(b)=-\frac{qh}{b^2+q^2+bgq}
\end{equation}
because $h^2+g^2/4=1$. In a similar way, if 
$$h(b):=\arctan \frac{2b+gq}{2hq}$$ 
then
\begin{equation}
\label{h}
\lim_{b\to \pm \infty}h(b)=\pm \frac{\pi}{2}\ \ \textrm{and}\ \ h'(b)=\frac{qh}{b^2+q^2+bgq}.
\end{equation}
Therefore $f+h$ is constant on the connected parts of the domain. Taking the limits $b\to \infty$ and $b\to -\infty$, respectively, we have 
$$\arctan \frac{q+\frac{g}{2}b}{hb}+\arctan \frac{2b+gq}{2hq}=\left\{\begin{array}{rl}
&+\frac{\pi}{2}+\arctan \frac{G}{2}\ \ \textrm{if} \ \ b>0\\
&\\
&-\frac{\pi}{2}+\arctan \frac{G}{2}\ \ \textrm{if} \ \ b < 0.
\end{array}
\right.$$
Therefore
\begin{equation}
\label{metrics}
\Phi=\arctan \frac{2b+gq}{2hq}.
\end{equation}
The following table shows the correspondence between the notations.

\vspace{0.5cm}
\begin{center}
\begin{tabular}{|c|c|}
\hline
&\\
Finsleroid-Finsler metric \cite{A2} & formula (\ref{main})\\
&\\
\hline
&\\
$b=b_i({\bf x})y^i$ (the Finsleroid axis $1$ - form) & $d\alpha$\\
&\\
\hline
&\\
$a_{ij}({\bf x})y^iy^j$ (the Riemannian norm sqare)& $F_*^2$\\
&\\
\hline
&\\
$g=g({\bf x})$ and $-2 < g({\bf x})< 2$ (the Finsleroid charge)& $\frac{K({\bf x})}{2}$\\
&\\
\hline
&\\
$h=\sqrt{1-\frac{\ g^2}{4}}$&$\sqrt{1-\frac{\ K^2}{16}}=\frac{\sqrt{16-K^2}}{4}$\\
&\\
\hline
&\\
$G=g/h$ & $\frac{2K}{\sqrt{16-K^2}}$\\
&\\
\hline
\end{tabular}
\end{center}

\vspace{0.5cm}
\noindent
We have that
$$\left(b^2+gqb+q^2\right)(v+tX_*(p))=t^2+\frac{K(p)}{2}F_*(v)t+F_*^2(v),$$
$$\Phi(v+tX_*(p))\stackrel{(\ref{metrics})}{=}\arctan \frac{2b+gq}{2hq}(v+tX_*(p))=\arctan \frac{1}{\sqrt{16-K^2(p)}}\left(\frac{4t}{F_*(v)}+K(p)\right)$$
and, consequently, 
\begin{equation}
\label{metrics1}
E(v+tX_*(p))\stackrel{(\ref{main})}{=}K^*(p)\left(F_*^2(v)+K(p)F_*(v)\frac{t}{2}+t^2\right)e^{2A(v,t)}=
\end{equation}
$$2K^*(p)e^{-G(p)\arctan \frac{G(p)}{2}}\cdot \ \textrm{Finsleroid-Finsler energy},$$
where
$$G=\frac{2K}{\sqrt{16-K^2}}$$
and the Finsleroid-Finsler energy is 
$$\frac{1}{2}e^{G\Phi}(b^2+gqb+q^2)$$
in the sense of formula (\ref{asanovmetric1}).
\end{Pf}

\section{The proofs of the main theorems}

\subsection{A continuity argument and the proof of Theorem B} In the previous section we proved that metric (\ref{main}) is conformal to a Finsleroid-Finsler metric. Since Finsleroid-Finsler metrics have indicatrices of constant positive curvature \cite{A2} (formula (2.32), page 284) it follows that the condition of the positivity of the sectional curvature in the main theorem can be extended to a local neighbourhood $U$ of $p$ together with the regularity condition $d_q \alpha\neq 0$ ($q\in U$) and the steps of the proof can be repeated to have the special form (\ref{main}) of $E_{|TU}$.  

\subsection{The converse of Theorem B and the proof of Theorem A} Using formulas (A.5), (A.7) and (A.12) in \cite{A2} (page 295) a direct computation shows that 
$$EX^l\alpha_l=Eg^{kl}\alpha_k\alpha_l=\left\{\begin{array}{rl}
&\frac{1}{2}(b^2+q^2)\ \ \textrm{(Asanov's notations)} \\
&\\
&\frac{1}{2}F_*^2=E_*.
\end{array}
\right.$$
This means that 
$$G^l\alpha_l-\tilde{G}^l\alpha_l\stackrel{(\ref{eq:6})}{=}EX^l\alpha_l-(y^m\alpha_m)^2$$ 
is quadratic in the tangent spaces, i.e.
$$0=\frac{\partial^3 G^l}{\partial y^i \partial y^j \partial y^k}\alpha_l-\frac{\partial^3 \tilde{G}^l}{\partial y^i \partial y^j \partial y^k}\alpha_l \ \ \ \Rightarrow\ \ \ \tilde{P}_{ijk}^l\alpha_l=P_{ijk}^l\alpha_l.$$
On the other hand formula (\ref{eq:6}) shows that 
$$\tP_{ijk}^l=P_{ijk}^l$$
if and only if 
$$0=\frac{\partial^3 G^l}{\partial y^i \partial y^j \partial y^k}-\frac{\partial^3 \tilde{G}^l}{\partial y^i \partial y^j \partial y^k},$$
i.e. $EX^l$ is quadratic in the tangent spaces for all indices $l=1, \ldots, n$. Since
\begin{equation}
\label{initial2}
EX^l=Eg^{kl}\alpha_k=\frac{1}{2}(b^2+q^2)-\frac{gq}{2}v^l,
\end{equation}
where $v^l=y^l-bb^l$ (see Asanov's notations (A.5) in \cite{A2}), we have that it is impossible because of the term $q$ unless $g=0$ (the Finslerian charge vanishes) and the space is Riemannian. For different proofs of Theorem A and the solution of Matsumoto's problem see \cite{V3}, \cite{V6} and \cite{V7}.

\section{A note on the two-dimensional case}

The case of dimension $2$ seems to be easier in the beginning: Theorem 1 and Theorem 2 is automatically holds because three vector fields must be linearly dependent in a two-dimensional vector space (vector fields which are tangential to the indicatrix form a one-dimensional linear space; cf. Remark 2 and Theorem 1). The rigidity conditions of type
$$\tP_{ijk}^l=P_{ijk}^l\ \ \textrm{or}\ \ \tP_{ijk}^l\alpha_l=P_{ijk}^l\alpha_l$$
implies the existence of the associated Riemannian metric (see section 4.1.) independently of the dimension of the space and differential equation (\ref{diffeqversion2dricatti}) follows in the same way. Unfortunately, the solution seems to be hard because of two main reasons:
\begin{itemize}
\item we need at least three independent directions to conclude that $b^2=1$ in section 4.4,
\item  we should pay a special attention to the initial conditions in section 4.5: they must be formulated along the one-dimensional linear subspace
$$T_pN=\{ v\in T_pM\ | \ v(\alpha)=0 \},\ \ \textrm{where}\ \ \dim T_pM=1.$$
This means that the origin, as the singularity of the Finslerian setting, divides the subspace into two disjoint connected parts (open half lines) and the initial condition for the unknown function 
$$z(t)=\frac{X_*^l\partial E /\partial y^l}{E}(v+tX_*(p))$$
can be described independently along the opposite directions $\pm v$ (cf. Theorem 5 in dimension $3$). Using $b^2=1$ as {\bf an additional condition} we can linearize the term $gq$ on the right hand side of formula (\ref{initial2}) by choosing the "constant" with opposite signs $\pm g$ for the opposite half-lines. Therefore singular solutions of the generalized Matsumoto's problem in dimension $2$ can be presented. Especially, any two-dimensional Finsler space with constant main scalar admits conformal changes  keeping the mixed curvature tensor of the Berwald connection invariant; for the details see \cite{VV} and Berwald's list of Finsler spaces with constant main scalar \cite{Berwald1}, formulas 118 I-III; see also \cite{Berwald2}.
\end{itemize}

\end{document}